\documentclass[12pt,thmsa]{article}%
\usepackage{amsmath}
\usepackage{amsfonts}
\usepackage{amssymb}
\usepackage{sw20bams}%
\usepackage{graphicx}
\providecommand{\U}[1]{\protect\rule{.1in}{.1in}}
\begin{document}

\author{Steven Finch}
\title{Least Capacity Point of Triangles }
\date{July 15, 2014}
\maketitle

\begin{abstract}
Let $D$ be a compact convex domain in the plane. \ P\'{o}lya \&\ Szeg\"{o}
and, independently, Levi \&\ Pan defined the point $p\in D$ that is
\textquotedblleft best insulated from the boundary $C$ of $D$%
\textquotedblright. \ We compute $p$ in the case when $C$ is an isosceles
right triangle, revisiting exact results from the study of complex conformal mappings.

\end{abstract}

\footnotetext{Copyright \copyright \ 2014 by Steven R. Finch. All rights
reserved.}Let $T=\left\{  x+iy\in\mathbb{C}:x>0\text{, }y>0\text{,
}x+y<1\right\}  $, the interior of a right isosceles triangle. \ Let $\Delta$
denote the unit disk in $\mathbb{C}$. \ Let $w\in T$. \ The Riemann mapping
theorem guarantees the existence of a conformal map $f_{w}:T\rightarrow\Delta$
such that $f_{w}(w)=0$. Moreover $f_{w}$ extends continuously to the boundary
of $T$, mapping it homeomorphically onto the unit circle.

P\'{o}lya \&\ Szeg\"{o} \cite{PS} wrote about the \textbf{inner radius}
$r_{w}$ of $T$ \textbf{relative to} $w$; the point $w$ that maximizes $r_{w}$
is the same as the \textbf{least capacity point} characterized sixty years
later by Levi \&\ Pan \cite{LP}. \ It is sufficient to examine the ratio
$f_{w}(z)/(z-w)$ or, more precisely, to minimize the difference%
\[
\underset{\text{Green's function for }T}{\underbrace{\ln\left\vert
f_{w}(z)\right\vert }}-\;\;\;\ln\left\vert z-w\right\vert
\]
in the limit as $z\rightarrow w$, over all $w\in T$. \ The first term is
harmonic in $T-\{w\}$, is bounded outside every neighborhood of $w$, and
vanishes at the boundary \cite{R1}; the additional term serves to repair the
singularity at $w$.

Our purpose is to find an explicit expression for $f_{w}$. \ We begin with a
beautiful formula communicated in \cite{R2}. \ This is followed by commentary
regarding specific examples in the literature. \ We conclude with a second
approach, drawing upon results in \cite{Ko, Mt}.

The quantity $r_{w}$ is the same as the \textbf{conformal radius} relative to
$w$ in our simple setting (involving triangles only).

Two addenda have emerged over time, one concerning the $30^{\circ}$%
-$60^{\circ}$-$90^{\circ}$ triangle and the other concerning the (fairly
arbitrary) $6$-$9$-$13$ triangle. \ We wonder whether a least capacity point
qualifies as a \textit{triangle center} \cite{Kb}. \ The present paper, in
part, continues our discussion \cite{Fi} in which three other candidate
triangle centers were examined. \ 

\section{Weierstrass Sigma}

Using the Schwarz reflection principle, $f_{w}$ extends across the hypotenuse
of $T$ to an analytic function (still called $f_{w}$) on the square $\left\{
x+iy\in\mathbb{C}:0<x<1\text{, }0<y<1\right\}  $ except for a simple pole at
$w^{\prime}=1+i-i\overline{w}$. \ Using the Schwarz reflection principle twice
more, $f_{w}$ extends to the doubled square $\left\{  x+iy\in\mathbb{C}%
:-1<x<1\text{, }-1<y<1\right\}  $ with zeroes at $\left\{  \pm w,\pm
\overline{w^{\prime}}\right\}  $ and poles at $\left\{  \pm w^{\prime}%
,\pm\overline{w}\right\}  $. Clearly $f_{w}(-1+iy)=f_{w}(1+iy)$ and
$f_{w}(x-i)=f_{w}(x+i)$ always. Let $\Lambda$ denote the lattice $\left\{
2m+2ni:m,n\in\mathbb{Z}\right\}  $. \ We may extend $f_{w}$ to a doubly
periodic meromorphic function on $\mathbb{C}$ via $f_{w}(z+\lambda)=f_{w}(z)$
for all $\lambda\in\Lambda$. \ Assume that $f_{w}(1)=1$ without loss of
generality. \ We deduce that \cite{Co}%
\[%
\begin{array}
[c]{ccc}%
f_{w}(z)=C_{w}\dfrac{\sigma(z-w)\sigma(z+w)\sigma(z-\overline{w^{\prime}%
})\sigma(z+\overline{w^{\prime}})}{\sigma(z-w^{\prime})\sigma(z+w^{\prime
})\sigma(z-\overline{w})\sigma(z+\overline{w})} &  & \text{for all }z\in T
\end{array}
\]
where
\[
\sigma(z)=z%
{\displaystyle\prod\limits_{\substack{\lambda\in\Lambda\\\lambda\neq0 }}}
\left(  1-\frac{z}{\lambda}\right)  \exp\left(  \frac{z}{\lambda}+\frac{z^{2}%
}{2\lambda^{2}}\right)  ,
\]%
\[
C_{w}=\frac{\sigma(1-w^{\prime})\sigma(1+w^{\prime})\sigma(1-\overline
{w})\sigma(1+\overline{w})}{\sigma(1-w)\sigma(1+w)\sigma(1-\overline
{w^{\prime}})\sigma(1+\overline{w^{\prime}})}.
\]
In order to employ Mathematica (or other computer algebra package), the
half-periods $1$, $i$ give rise to invariants%
\[%
\begin{array}
[c]{ccc}%
g_{2}=\dfrac{1}{256\pi^{2}}\Gamma\left(  \dfrac{1}{4}\right)  ^{8}%
=11.8170450080..., &  & g_{3}=0
\end{array}
\]
which must be passed to the software implementation of $\sigma$. As
$z\rightarrow w$, the ratio $\left\vert f_{w}(z)/(z-w)\right\vert $ simplifies
to%
\[
h(w)=\left\vert \frac{\sigma(1-w^{\prime})\sigma(1+w^{\prime})\sigma
(1-\overline{w})\sigma(1+\overline{w})}{\sigma(1-w)\sigma(1+w)\sigma
(1-\overline{w^{\prime}})\sigma(1+\overline{w^{\prime}})}\cdot\dfrac
{\sigma(2w)\sigma(w-\overline{w^{\prime}})\sigma(w+\overline{w^{\prime}}%
)}{\sigma(w-w^{\prime})\sigma(w+w^{\prime})\sigma(w-\overline{w}%
)\sigma(w+\overline{w})}\right\vert
\]
because $\sigma(z-w)/(z-w)\rightarrow1$. \ Further simplification does not
seem possible. \ Numerical minimization gives the least capacity point to be
$w_{0}=(1+i)t_{0}$, where%
\[
t_{0}=0.3011216108413220815538254....
\]
No closed-form expression for $t_{0}$ is apparent, at least not here. \ We
observe (to high precision) that%
\[
\frac{1}{h(w_{0})}=0.3346161009568417919464744...=\frac{4\sqrt{2\pi}}{3^{3/4}%
}\Gamma\left(  \dfrac{1}{4}\right)  ^{-2}%
\]
which is encouraging since the latter is the maximum inner radius for $T$
\cite{PS}. \ A\ rigorous proof, however, remains open.

\section{Weierstrass P}

The zero $w_{0}$ in the preceding section identified a conformal map
$f_{w_{0}}:T\rightarrow\Delta$ that is, in particular, extremal in some sense.
\ Any map onto the upper half plane $\mathbb{C}^{+}$ can be easily recast as a
map onto the disk $\Delta$ (via composition with a linear fractional
transformation). \ As a slight detour, let us similarly identify other
better-known conformal maps $T\rightarrow\mathbb{C}^{+}$ that have appeared in
the literature. \ 

Two maps $\varphi$, $\psi$ are prescribed to take the following values on the
vertices of $T$:%
\[
(0,1,i)\overset{\varphi}{\longmapsto}(0,1,\infty),
\]%
\[
(0,1,i)\overset{\psi}{\longmapsto}(\infty,0,1).
\]
The Schwarz-Christoffel transformation gives \cite{Sp, Ky, Nh}%
\[
\varphi^{(-1)}(\zeta)=\frac{1}{\sqrt{2\pi}}\Gamma\left(  \dfrac{1}{4}\right)
^{2}B\left(  \zeta,\frac{1}{2},\frac{1}{4}\right)  ,
\]%
\[
\psi^{(-1)}(\zeta)=\frac{i-1}{\sqrt{\pi}}\Gamma\left(  \dfrac{1}{4}\right)
^{2}B\left(  \zeta,\frac{1}{4},\frac{1}{4}\right)  +1
\]
where%
\[
B(\zeta,\alpha,\beta)=%
{\displaystyle\int\limits_{0}^{\zeta}}
s^{\alpha-1}(1-s)^{\beta-1}ds
\]
is the incomplete Euler beta function. \ Although expressions for
$\varphi^{(-1)}$, $\psi^{(-1)}$ are famous, their inverses are comparatively
obscure. \ Geyer \cite{Gy} provided%
\[
\psi(z)=-\frac{1}{4\wp(1)}\frac{(\wp(z)-\wp(1))^{2}}{\wp(z)}%
\]
where%
\[
\wp(z)=\frac{1}{z^{2}}+%
{\displaystyle\sum\limits_{\substack{\lambda\in\Lambda\\\lambda\neq0}}}
\left(  \frac{1}{(z-\lambda)^{2}}-\frac{1}{\lambda^{2}}\right)  ,
\]%
\[
\wp(1)=\frac{1}{32\pi}\Gamma\left(  \dfrac{1}{4}\right)  ^{4}=1.7187964545....
\]
Choose the linear fractional transformation $\mathbb{C}^{+}\rightarrow\Delta$
to be $\zeta\longmapsto(\zeta-i)/(\zeta+i)$. \ Hence we wish to solve the
equation $\psi(w)=i$, but this is immediately seen to yield
\[
w=\frac{i-1}{\sqrt{\pi}}\Gamma\left(  \dfrac{1}{4}\right)  ^{2}B\left(
i,\frac{1}{4},\frac{1}{4}\right)  +1=0.2970894700...+(0.1926647354...)i.
\]
It is clear that%
\[
\overline{\psi(i\bar{z})}\cdot\varphi(z)=1
\]
and we wish to solve the equation $\varphi(w)=i$, but this too is immediately
seen to yield%
\[
w=\frac{1}{\sqrt{2\pi}}\Gamma\left(  \dfrac{1}{4}\right)  ^{2}B\left(
i,\frac{1}{2},\frac{1}{4}\right)  =0.1926647354...+(0.2970894700...)i.
\]
A representation of arbitrary $f_{w}$ in terms of $\wp$ (analogous to our
representation in terms of $\sigma$) is also possible \cite{Co}, but we
haven't pursued this. \ 

\section{Jacobi Elliptic}

Let
\[
F[\phi,m]=%
{\displaystyle\int\limits_{0}^{\sin(\phi)}}
\dfrac{d\tau}{\sqrt{1-\tau^{2}}\,\sqrt{1-m\,\tau^{2}}}\,
\]
denote the incomplete elliptic integral of the first kind and $K[m]=F[\pi
/2,m]$; we purposefully choose formulas here to be consistent with
Mathematica. The three basic Jacobi elliptic functions are defined via
\begin{align*}
u  & =%
{\displaystyle\int\limits_{0}^{\operatorname*{sn}(u,m)}}
\frac{d\tau}{\sqrt{1-\tau^{2}}\,\sqrt{1-m\,\tau^{2}}}=%
{\displaystyle\int\limits_{\operatorname*{cn}(u,m)}^{1}}
\frac{d\tau}{\sqrt{1-\,\tau^{2}}\,\sqrt{m\,\tau^{2}+(1-m)}}\\
& =%
{\displaystyle\int\limits_{\operatorname*{dn}(u,m)}^{1}}
\frac{d\tau}{\sqrt{1-\,\tau^{2}}\,\sqrt{\tau^{2}-(1-m)}}%
\end{align*}
and we shall require all three of these. \ Define \
\[
\kappa=\frac{K[1/2]}{\sqrt{2}}=\frac{1}{2^{5/2}\sqrt{\pi}}\Gamma\left(
\dfrac{1}{4}\right)  ^{2}=1.3110287771...=\frac{1.8540746773...}{\sqrt{2}}%
\]
and let $\tilde{T}=\left\{  x+iy\in\mathbb{C}:y>0\text{, }y<x+\kappa\text{,
}y<-x+\kappa\right\}  $. \ As in Section 2, let us first examine a conformal
map\ $\theta:\tilde{T}\rightarrow\mathbb{C}^{+}$ which takes prescribed values
on the vertices of $\tilde{T}$:
\[
(-\kappa,\kappa,i\kappa)\overset{\theta}{\longmapsto}(-1,1,\infty).
\]
The Schwarz-Christoffel transformation gives \cite{NP}%
\[
\theta^{(-1)}(\zeta)=\frac{1}{2}%
{\displaystyle\int\limits_{0}^{\zeta}}
\frac{ds}{\left(  1-s^{2}\right)  ^{3/4}}=\frac{\zeta}{2}\,_{2}F_{1}\left(
\frac{1}{2},\frac{3}{4},\frac{3}{2},\zeta^{2}\right)
\]
which involves the following Gauss hypergeometric function:
\[
_{2}F_{1}\left(  \frac{1}{2},\frac{3}{4},\frac{3}{2},z\right)  =\frac
{\Gamma(1/4)}{2\sqrt{2}\pi}%
{\displaystyle\sum\limits_{n=0}^{\infty}}
\frac{\Gamma(n+1/2)\Gamma(n+3/4)}{\Gamma(n+3/2)}\frac{z^{n}}{n!}.
\]
Although the expression for $\theta^{(-1)}$ is famous, its inverse is
comparatively obscure. \ Kober \cite{Ko, Mt} provided%
\[
\theta(z)=\sqrt{2}\operatorname*{sn}\left(  \sqrt{2}z,\frac{1}{2}\right)
\operatorname*{dn}\left(  \sqrt{2}z,\frac{1}{2}\right)  .
\]
Choose the linear fractional transformation $\mathbb{C}^{+}\rightarrow\Delta$
to be $\zeta\longmapsto(\zeta-i)/(\zeta+i)$. \ Thus we wish to solve the
equation $\theta(w)=i$, but this is immediately seen to yield
\[
w=\frac{i}{2}\,_{2}F_{1}\left(  \frac{1}{2},\frac{3}{4},\frac{3}{2},-1\right)
=(0.4154481080...)i=(0.3168871006...)\kappa i.
\]

Now let $w\in\tilde{T}$ be arbitrary. \ Define\ $f_{w}:\tilde{T}%
\rightarrow\Delta$ for which $f_{w}(w)=0$ via \cite{Nh}%
\begin{align*}
f_{w}(z)  & =-\frac{\frac{\theta(z)-i}{\theta(z)+i}-\frac{\theta(w)-i}%
{\theta(w)+i}}{1-\left(  \frac{\theta(z)-i}{\theta(z)+i}\right)
\overline{\left(  \frac{\theta(w)-i}{\theta(w)+i}\right)  }}\\
& =\frac{\theta(z)-\theta(w)}{\theta(z)-\overline{\theta(w)}}%
\end{align*}
(following the construction of Green's function in \cite{KS}, but beware of
misprints). \ As $z\rightarrow w$, the ratio $\left\vert f_{w}%
(z)/(z-w)\right\vert $ simplifies to%
\[
h(w)=\left\vert \frac{\theta^{\prime}(w)}{\theta(w)-\overline{\theta(w)}%
}\right\vert
\]
where $\theta^{\prime}$ denotes the derivative of $\theta$. \ Numerical
minimization gives the least capacity point to be%
\begin{align*}
\tilde{w}_{0}  & =(0.3977567783173558368923490...)\kappa i\\
& =(1-2t_{0})\kappa i
\end{align*}
as expected, since
\[
\frac{\kappa-\left\vert \tilde{w}_{0}\right\vert }{\sqrt{2}\kappa}=\frac
{\sqrt{2}t_{0}}{1}%
\]
by the similarity of triangles $\tilde{T}$ and $T$. \ Restricting attention to
the $y$-axis only, we have \
\[
h(iy)=\frac{i\operatorname*{cn}\left(  i\sqrt{2}y,\frac{1}{2}\right)  ^{3}%
}{\sqrt{2}\operatorname*{sn}\left(  i\sqrt{2}y,\frac{1}{2}\right)
\operatorname*{dn}\left(  i\sqrt{2}y,\frac{1}{2}\right)  }.
\]
Differentiating with respect to $y$ and setting the result equal to zero, the
equation%
\[
\operatorname*{dn}\left(  i\sqrt{2}y,\frac{1}{2}\right)  =\sqrt{\frac
{1+\sqrt{3}}{2}}%
\]
is found, therefore%
\[
t_{0}=\operatorname{Re}\left\{  \frac{1}{2\kappa}F\left[  \arcsin\sqrt
{\frac{1+\sqrt{3}}{2}},2\right]  \right\}  =0.3011216108413220815538254...
\]
is the sought-after closed-form expression. \ Such an outcome was not apparent
in Section 1.

An old paper by Love \cite{Lv} discusses conformal maps on four exceptional
triangles (including the isosceles right triangle) and utilizes the $\wp$
function; unfortunately we haven't succeeded in following the details. Two
other papers\ \cite{CGS, SRS}, despite promising titles, evidently assess
Green's function created for different settings (non-Laplacian) than ours.

\section{Addendum:\ $30^{\circ}$-$60^{\circ}$-$90^{\circ}$ Triangle}

Without entering any lengthy explanations, let \
\[
\kappa=\frac{1}{2^{5/3}\sqrt{\pi}}\Gamma\left(  \dfrac{1}{3}\right)
\Gamma\left(  \dfrac{1}{6}\right)  =\frac{5.2999162508...}{2}%
\]
and define $T=\left\{  x+iy\in\mathbb{C}:x>0\text{, }y>0\text{, }\sqrt
{3}x+y<\sqrt{3}\kappa\right\}  $. \ The conformal map\ $\theta:T\rightarrow
\mathbb{C}^{+}$ taking prescribed values on the vertices of $T$:
\[
(0,\kappa,i\sqrt{3}\kappa)\overset{\theta}{\longmapsto}(0,1/4,\infty)
\]
is given by
\[
\theta(z)=\frac{3\sqrt{3}\operatorname*{sn}\left(  \dfrac{2^{2/3}}{3^{3/4}%
}z,\dfrac{2+\sqrt{3}}{4}\right)  ^{2}\operatorname*{dn}\left(  \dfrac{2^{2/3}%
}{3^{3/4}}z,\dfrac{2+\sqrt{3}}{4}\right)  ^{2}}{\left\{  1+\operatorname*{cn}%
\left(  \dfrac{2^{2/3}}{3^{3/4}}z,\dfrac{2+\sqrt{3}}{4}\right)  \right\}
^{4}}.
\]
(We could not decipher pp. 186--187 of \cite{Ko}, hence we turned to pp.
184--185 and made necessary adjustments.) \ For example, $\theta(w)=i$ occurs
when
\begin{align*}
w  & =3\left(  \frac{1}{2}+\frac{1+i}{\sqrt{2}}\right)  ^{1/3}\,_{2}%
F_{1}\left(  \frac{1}{3},\frac{2}{3},\frac{4}{3},\frac{1}{2}+\frac{1+i}%
{\sqrt{2}}\right)  -\kappa\\
& =0.7065812599...+(1.6814450943...)i\\
& =(0.2666386510...)\kappa+(0.6345176092...)\kappa i.
\end{align*}

Now let $w\in T$ be arbitrary. \ Define\ $f_{w}:T\rightarrow\Delta$ for which
$f_{w}(w)=0$ as before; the ratio $\left\vert f_{w}(z)/(z-w)\right\vert $
tends to%
\[
h(w)=\left\vert \frac{\theta^{\prime}(w)}{\theta(w)-\overline{\theta(w)}%
}\right\vert
\]
as $z\rightarrow w$. \ Numerical minimization gives the least capacity point
to be%
\begin{align*}
w_{0}  & =(0.3599371272406945147550792...)\kappa+\\
& (0.4062604057445303763104149...)\kappa i.
\end{align*}
We have not attempted to find a closed-form expression for $w_{0}$. \ To high
precision,
\[
\frac{1}{h(w_{0})}=(0.2105704622445114724079460...)(2\kappa)=\frac{2^{4/3}\pi
}{5^{5/12}}\Gamma\left(  \dfrac{1}{3}\right)  ^{-3}(2\kappa)
\]
which is the (corrected) maximum inner radius for $T$ \cite{PS}.

\section{Addendum:\ $6$-$9$-$13$ Triangle}

No such exact formulas can be found for arbitrary triangles. \ The
Schwarz-Christoffel toolbox for Matlab \cite{DT, Dr}, coupled with the
Optimization toolbox, makes numerical computations of least capacity points
readily accessible. \ Recall that we wish to assess whether such points can be
treated as triangle centers. \ For the triangle with vertices%
\[
0,\;\;\;6,\;\;\;-\frac{13}{3}+\frac{4\sqrt{35}}{3}i
\]
the following code:%

\[%
\begin{array}
[c]{l}%
\text{\texttt{function q = arbitra(w)}}\\
\text{\texttt{p = polygon([0 6 -13/3+(4*sqrt(35)/3)*i])}}\\
\text{\texttt{f = diskmap(p,scmapopt('Tolerance',1e-18))}}\\
\text{\texttt{f = center(f,w(1)+i*w(2));}}\\
\text{\texttt{p = parameters(f);}}\\
\text{\texttt{q = -abs(p.constant);}}%
\end{array}
\]
gives (for example) that the inner radius at the centroid is%
\[
-\operatorname*{arbitra}\left(  \frac{5}{9}+\frac{4\sqrt{35}}{9}i\right)
=1.802305....
\]
Using the centroid as a starting guess, we solve a constrained minimization
problem as follows:
\[%
\begin{array}
[c]{l}%
\text{\texttt{format long}}\\
\text{\texttt{options=optimset('Algorithm','interior-point','TolCon',
1e-15);}}\\
\text{\texttt{A = [-4*sqrt(35) -13; 0 -1; 4*sqrt(35) 31]}}\\
\text{\texttt{b = [0 0 24*sqrt(35)]}}\\
\text{\texttt{v0 = [5/9 4*sqrt(35)/9];}}\\
\text{\texttt{[v,fv] = fmincon(@arbitra,v0,A,b,[],[],[],[],[],options);}}%
\end{array}
\]
yielding the maximum inner radius to be $1.979479...$ and the corresponding
least capacity point to be $0.929617...+(1.842564...)i$. \ This particular
triangle serves as a benchmark in \cite{Kb} to distinguish various centers.
\ The imaginary part is the perpendicular distance from the proposed center to
the shortest triangle side. \ Since the numerical value $1.842...$ does not
appear in the database, we infer that this center is new.

Figures 1, 2, 3 provide conformal map images of ten evenly-spaced concentric
circles in the disk. \ These are optimal in the sense that their center is
\textquotedblleft best insulated\textquotedblright\ from the triangle
boundary. Orthogonal trajectories are also indicated. \ 

The literature on this subject is larger than we originally thought. \ The
phrase \textbf{conformal center} is sometimes used to denote what we call the
least capacity point. \ (This is not to be confused with a different sense of
the same phrase in \cite{DT, Dr}.) \ Some discussion of relevant numerical
optimization based on the Schwarz-Christoffel transformation occurred years
ago \cite{Fl}. \ Precise inequalities relating radii and various points have
also been formulated \cite{Kz}.

The same phrase is used to denote yet another triangle center in \cite{Ia}.
\ Starting from such a location, a particle undergoing Brownian motion is
equally likely to exit through any of the triangle sides. \ As far as is
known, this topic is distinct from our study. \ Certain integrals and series
in \cite{Ia} deserve greater attention.%

\begin{figure}[ptb]%
\centering
\includegraphics[
height=4.2704in,
width=4.4512in
]%
{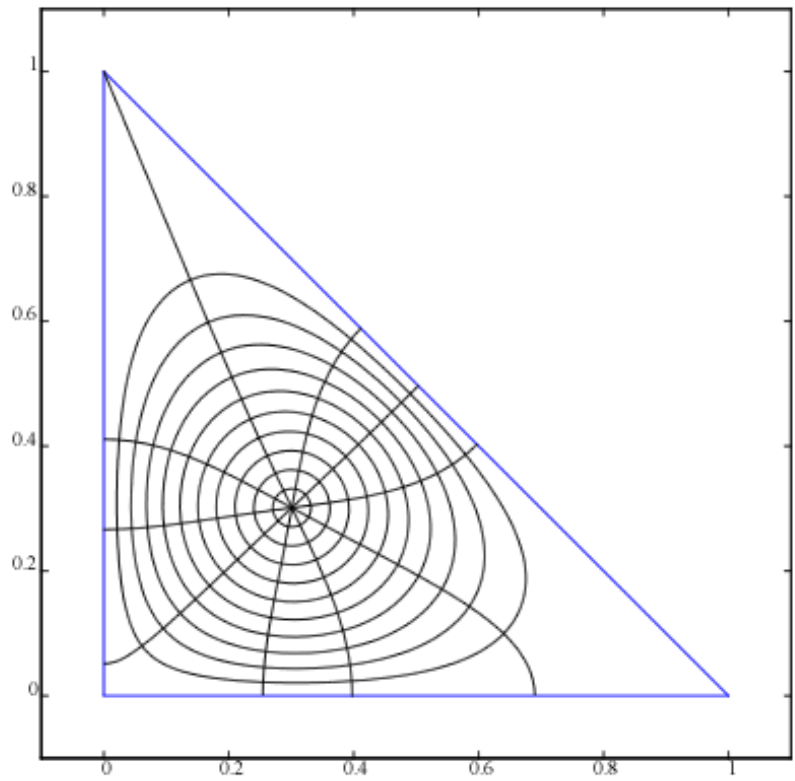}%
\caption{Images of ten concentric circles, center at $0.301+(0.301)i$.}%
\end{figure}
%

\begin{figure}[ptb]%
\centering
\includegraphics[
height=4.2704in,
width=4.3145in
]%
{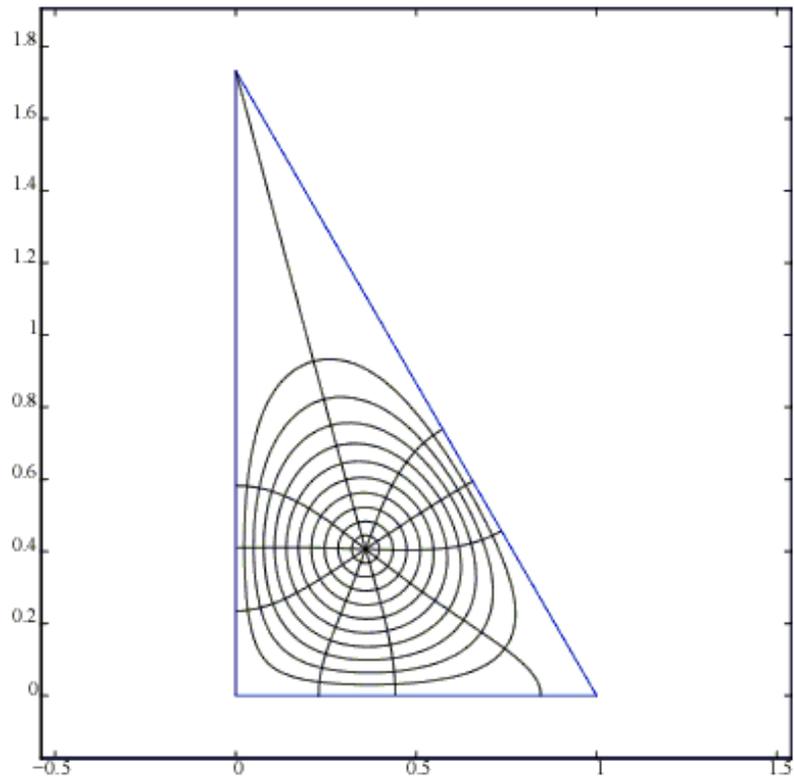}%
\caption{Images of ten concentric circles, center at $0.359+(0.406)i$.}%
\end{figure}
\ \
\begin{figure}[ptb]%
\centering
\includegraphics[
height=4.2696in,
width=4.3379in
]%
{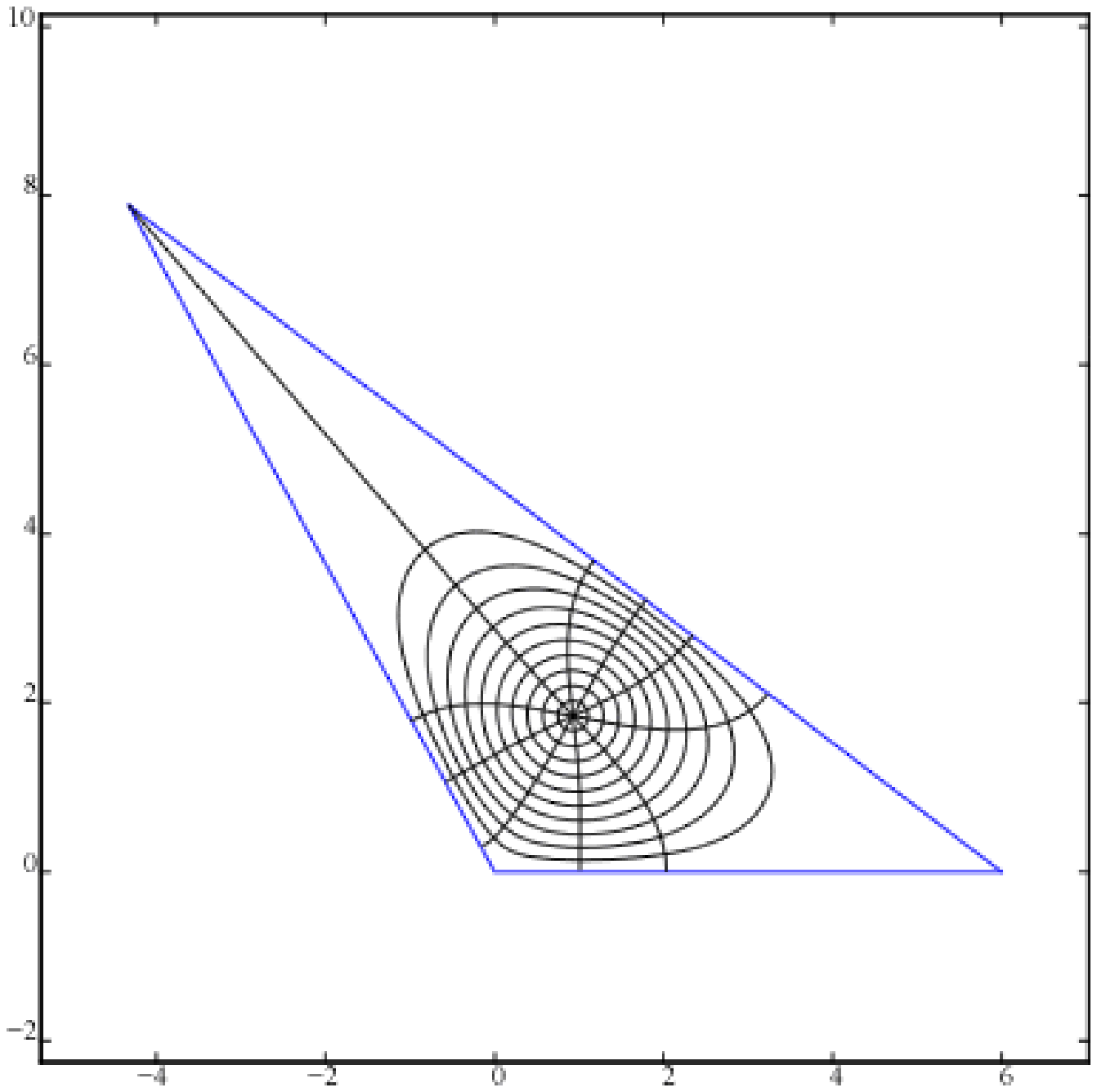}%
\caption{Images of ten concentric circles, center at $0.929+(1.842)i$.\ }%
\end{figure}

\section{Acknowledgements}

I am grateful to Thomas Ransford \cite{R1, R2} for providing the expression
for $f_{w}$ involving the Weierstrass sigma function, along with detailed
proofs of theorems and answers to several questions.


\begin{thebibliography}{99}                                                                                               %
\bibitem {PS}G. P\'{o}lya and\ G. Szeg\"{o}, \textit{Isoperimetric
Inequalities in Mathematical Physics}, Princeton Univ. Press, 1951, pp. 1--3,
256--258; MR0043486 (13,270d).

\bibitem {LP}M. Levi and J. Pan, Minimal capacity points and the lowest
eigenfunctions, http://arxiv.org/abs/1104.0555.

\bibitem {R1}T. Ransford, \textit{Potential Theory in the Complex Plane},
Cambridge Univ. Press, 1995, pp. 106--116; MR1334766 (96e:31001).

\bibitem {R2}T. Ransford, Green's function for $T$, unpublished note (2014).

\bibitem {Ko}H. Kober, \textit{Dictionary of Conformal Representations},
Dover, 1952, pp. 182--187; MR0049326 (14,156d).

\bibitem {Mt}G. Moretti, \textit{Functions of a Complex Variable},
Prentice-Hall, 1964, pp. 359--361, 375--376, 385--387.

\bibitem {Kb}C. Kimberling, Encyclopedia of Triangle Centers, http://faculty.evansville.edu/ck6/encyclopedia/.

\bibitem {Fi}S. Finch, In limbo: Three triangle centers, http://arxiv.org/abs/1406.0836.

\bibitem {Co}E. T. Copson, \textit{An Introduction to the Theory of Functions
of a Complex Variable}, Oxford Univ. Press, 1935, pp. 365--366, 371--372.

\bibitem {Sp}M. R. Spiegel, \textit{Complex Variables}, McGraw-Hill, 1964, pp.
204, 209, 221--222.

\bibitem {Ky}P. K. Kythe, \textit{Computational Conformal Mapping},
Birkh\"{a}user Boston, 1998, pp. 54--55, 76--77; MR1651941 (99k:65027).

\bibitem {Nh}Z. Nehari,\textit{\ Conformal Mapping}, Dover, 1975, pp. 182,
194--195; MR0377031 (51 \#13206).

\bibitem {Gy}L. Geyer, Conformal mapping from triangle to upper half plane in
terms of Weierstrass, http://math.stackexchange.com/questions/246309/conformal-mapping-from-triangle-to-upper-half-plane-in-terms-of-weierstrass-wp.

\bibitem {NP}R. Nevanlinna and V. Paatero, \textit{Introduction to Complex
Analysis}, Chelsea, 1982, pp. 332--335; MR0657146 (83d:30002).

\bibitem {KS}N. Kurt and M. Sezer, Solution of Dirichlet problem for a
triangle region in terms of elliptic functions, \textit{Appl. Math. Comput.}
182 (2006) 73--81; MR2292020 (2007k:35073).

\bibitem {Lv}A. E. H. Love, Vortex motion in certain triangles, \textit{Amer.
J. Math.} 11 (1889) 158--171.

\bibitem {CGS}J. Cuenca, F. Gautier and L. Simon, The image source method for
calculating the vibrations of simply supported convex polygonal plates,
\textit{J. Sound Vibration} 322 (2009) 1048--1069.

\bibitem {SRS}A. Scalia, A. Rigano and M. A. Sumbatyan, Reconstruction of
voids in elastic isosceles right triangle, \textit{Mechanics Research
Communications} 37 (2010) 650--654.

\bibitem {DT}T. A. Driscoll and L. N. Trefethen, \textit{Schwarz-Christoffel
Mapping}, Cambridge Univ. Press, 2002, pp. 9--18; MR1908657 (2003e:30012).

\bibitem {Dr}T. A. Driscoll, Schwarz-Christoffel Toolbox for Matlab,
http://www.math.udel.edu/\symbol{126}driscoll/SC/.

\bibitem {Fl}M. Flucher, An asymptotic formula for the minimal capacity among
sets of equal area, \textit{Calc. Var. Partial Differential Equations} 1
(1993) 71--86; MR1261718 (95c:35276).

\bibitem {Kz}V. O. Kuznetsov, On properties of the conformal radius of a
domain (in Russian), \textit{Zap. Nauchn. Sem. S.-Peterburg. Otdel. Mat. Inst.
Steklov. (POMI)} 276 (2001), \textit{Anal. Teor. Chisel i Teor. Funkts.} 17,
237--252, 352; Engl. transl. in \textit{J. Math. Sci. (N. Y.)} 118 (2003)
4871--4879; MR1850370 (2002h:30008).

\bibitem {Ia}A. Iannaccone, \textit{The Conformal Center of a Triangle or a
Quadrilateral}, BS\ thesis, Harvey Mudd College, 2003, https://www.math.hmc.edu/seniorthesis/archives/2003/aiannacc/.%

\begin{tabular}
[c]{lll}
& Steven Finch & \\
& Dept. of Statistics & \\
& Harvard University & \\
& Cambridge, MA, USA & \\
& \textit{steven\_finch@harvard.edu} &
\end{tabular}

\end{thebibliography}
\end{document}